\documentclass[]{IEEEtran4PSCC}
\ifCLASSINFOpdf
   \usepackage[pdftex]{graphicx}
\else
   \usepackage[dvips]{graphicx}
\fi

\usepackage[cmex10]{amsmath}
\usepackage{kpfonts}

\usepackage{amsmath}
\usepackage{mathtools}
\usepackage{amssymb}
\usepackage{amsthm}
\usepackage{amsfonts}
\usepackage{amscd}
\usepackage{times}
\usepackage{bm}
\usepackage{wrapfig}
\usepackage{pdfpages}
\usepackage{subfigure}
\usepackage{comment}
\usepackage{float}
\floatstyle{ruled}
\usepackage{textcomp}
\usepackage{wrapfig}
\usepackage{hyperref}
\hypersetup{
	colorlinks=true, 
	breaklinks=true,
	urlcolor= blue, 
	linkcolor= blue, 
	citecolor=red, 
	pdftitle={},
	pdfauthor={},
}
\usepackage{lipsum}
\usepackage{mwe}
\usepackage{amsfonts,amssymb,amsmath,textcomp,mathabx,empheq,amsthm}
\usepackage{algorithmicx}
\usepackage{algorithm}
\usepackage{algpseudocode}
\usepackage{accents}
\newfloat{model}{thp}{lop}
\floatname{model}{Model}
\newcommand\ubar[1]{\underaccent{\bar}{#1}}

\makeatletter
\let\old@ps@headings\ps@headings
\let\old@ps@IEEEtitlepagestyle\ps@IEEEtitlepagestyle
\def\psccfooter#1{%
    \def\ps@headings{%
        \old@ps@headings%
        \def\@oddfoot{\strut\hfill#1\hfill\strut}%
        \def\@evenfoot{\strut\hfill#1\hfill\strut}%
    }%
    \def\ps@IEEEtitlepagestyle{%
        \old@ps@IEEEtitlepagestyle%
        \def\@oddfoot{\strut\hfill#1\hfill\strut}%
        \def\@evenfoot{\strut\hfill#1\hfill\strut}%
    }%
    \ps@headings%
}
\makeatother

\psccfooter{%
        \parbox{\textwidth}{\hrulefill \\ \small{21st Power Systems Computation Conference} \hfill \begin{minipage}{0.2\textwidth}\centering \vspace*{4pt} \includegraphics[scale=0.06]{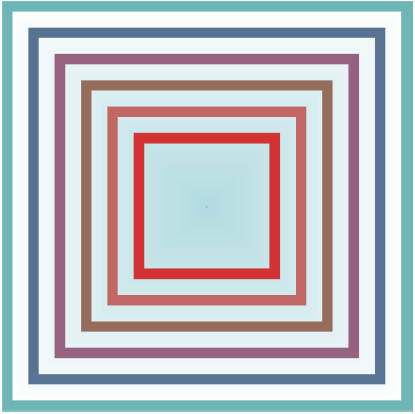}\\\small{PSCC 2020} \end{minipage} \hfill \small{Porto, Portugal --- June 29 -- July 3, 2020}}%
}
\begin{document}
\title{Proving Global Optimality of ACOPF Solutions}

\author{
S. Gopinath$^\dag$, H.L. Hijazi$^\dag$, T. Weisser$^\dag$, H. Nagarajan$^\dag$, M. Yetkin$^*$, K. Sundar$^\dag$, R.W. Bent$^\dag$ \\ 
$^\dag$Los Alamos National Laboratory, Los Alamos, NM, U.S.A. (Contact: hlh@lanl.gov) \\
$^*$Department of Industrial and Systems Engineering, Lehigh University, Bethlehem, PA, U.S.A.
}

\maketitle

\begin{abstract}
We present our latest contributions in terms of mathematical modeling and algorithm development for the global optimization of the ACOPF problem. These contributions allow us to close the optimality gap on a number of open instances in the PGLIB and NESTA benchmark libraries. This is achieved by combining valid cut generation with semidefinite programming-based bound tightening. The mathematical formulations along with the solution algorithms are implemented in the modeling framework Gravity (www.gravityopt.com), an open-source platform for reproducible numerical experiments.

\end{abstract}

\begin{IEEEkeywords}
ACOPF, Convex Relaxation, Global Optimization.
\end{IEEEkeywords}

\section{Introduction}

The AC Optimal Power Flow (ACOPF) problem lies at the heart of power systems computation. In recent years, a number of convex relaxations have been proposed with the intention of solving this NP-hard problem to global optimality \cite{qc_mpc_17,qc}.
These efforts led to efficient formulations and solution algorithms that can close the optimality gap on 95\% of instances found in the literature. In this work, we focus on the remaining 5\%.
We are particularly interested in pushing the reach of methods that do not require domain partitioning and spatial branching on continuous variables. We are concurrently pursuing other research directions based on our recent work \cite{lu2018tight} that rely on adaptive domain partitioning but this will not be covered here.


Our contributions are summarized below:
\begin{itemize}
    \item A new semidefinite programming-based bound tightening algorithm using determinant constraints \cite{polysdp}.
    \item Valid cut generation based on the Reformulation-Linearization Technique \cite{RLT}.
    \item The first open-source parallel multi-machine and multi-threaded implementation of the contributions above.
\end{itemize}


In Section \ref{sec:ACOPF}, we introduce the extended ACOPF problem formulation including line charging, bus shunts and transformers. 
In Section \ref{sec:SDP} we discuss Semidefinite Programming (SDP) relaxations of the ACOPF problem along with valid cut generation.
In Section \ref{sec:SDP-BT}, we introduce the new SDP-based bound tightening (SDP-BT) algorithm and in the last section we present our numerical experiments.

\section{ACOPF}
\label{sec:ACOPF}

\subsection{Notations} 
\label{sec:notation}
\vspace{10pt}
\resizebox{0.5\textwidth}{!}{%
\begin{tabular}{ll}
{\bf Grid Parameters:}\\
\vspace{0.2cm}
$\bm S^g_i = (\bm p^g_i, \bm q^g_i)$  	&Complex, active and reactive power generation at node $i$\\
$\bm S^d_{i} = (\bm p^d_i, \bm q^d_i)$  	&Complex, active and reactive power demand at node $i$ \\
$\bm c_{0i}, \bm c_{1i}, \bm c_{2i}$  			&Generation cost coefficients at node $i$ \\
$\bm t_{ij}$  				&Thermal limit along line $(i,j)$ \\
$\bm Y_{ij}$  				& Complex admittance of line $(i,j)$ $\bm Y_{ij}=\bm g_{ij}+ \bm i \bm b_{ij}$ \\
$\bm Y^c_{ij}$  				&Line charging admittance of line $(i,j)$  \\
$\bm Y^s_{i}$  				&Bus shunt admittance of bus $i$  \\
$\ubar{\bm S}^g_i, \bar{\bm S}^g_i$  	&$\ubar{\bm p}^g_i+\bm i \ubar{\bm q}^g_i, \bar{\bm p}^g_i+\bm i \bar{\bm q}^g_i $ \\
$\ubar{\bm v}_i, \bar{\bm v}_i$  	&Voltage magnitude bounds at node $i$\\
$\ubar{\bm \theta}_{ij} \ge -\pi/2, \bar{\bm \theta}_{ij} \le \pi/2$  &Phase angle difference bounds along line $(i,j)$\vspace{0.2cm}\\
{\bf Grid Variables:}\\
\vspace{0.2cm}
$S_{ij}= ( p_{ij}, q_{ij})$  	&Complex, active and reactive power flow from $i$ to $j$ \\
$V_{i}$  				&Complex voltage at node $i$\\
$I_{ij}$  				&Complex current flowing from $i$ to $j$ at end $i$\\
$l_{ij}$  				&Squared magnitude of current flowing from $i$ to $j$ at end $i$\\
$W_{ij}$    				& Lifted voltage product $V_iV_j^*$ \\
{\bf Other notations:}\\
\vspace{0.2cm}
$\bm {\mathcal G} = (N.E)$ 				&Graph with nodes $N$ and edges $E$\\
$\mathcal{R}(.)$     		& Real component of complex number\\
$\mathcal{I}(.)$     		& Imaginary component of complex number\\
${V}^*$     		& Conjugate of complex number\\
$|V|^2$     		& Square Magnitude of complex number $|V|^2 = \mathcal{R}(V)^2 + \mathcal{I}(V)^2$\\
$\ubar{\bm x}$ , $\bar{\bm x}$       &Lower and upper bounds on $x$\\
\end{tabular}
}

\vspace{11pt}
\subsection{Formulation} \label{sec:formulation}
The extended ACOPF problem is presented in Model \ref{model:ext_ac_opf}.
Note that complex inequalities correspond to component-wise constraints on the respective real and imaginary parts.

 \begin{model}[!h]
	\caption{Extended ACOPF}
	\label{model:ext_ac_opf}
	{\scriptsize
	\begin{subequations}
		\vspace{-0.2cm}
		\begin{align}
		& \mbox{\bf minimize: } \nonumber \\ 
		& \sum\limits_{i \in N} \bm c_{0i} + \bm c_{1i}\mathcal{R}(S^g_i) + \left(\bm c_{2i}\mathcal{R}(S^g_i)\right)^2 \nonumber\\
		& \mbox{\bf subject to: } \nonumber  \\
		& S_{ij} = (\bm Y^*_{ij}+ \bm Y^{\bm c^*}_{ij})\frac{| V_i |^2}{| \bm T_{ij}|^2}-\bm Y^*_{ij}\frac{V_iV^*_j}{\bm T_{ij}}  ,\quad \forall (i,j) \in E \label{eq:expsij} \\
			& S_{ji} = (\bm Y^*_{ij}+\bm Y^{c^*}_{ij})| V_j |^2-\bm Y^*_{ij}\frac{V_i^*V_j}{\bm T^*_{ij}} ,\quad \forall (i,j) \in E \label{eq:expsji}\\
            & S_i^g - \bm{S}_i^d -\bm{Y^s}_i| V_i |^2= \displaystyle\sum_{(i,j),(j,i) \in E} S_{ij} ,\quad \forall i \in N, \label{eq:exkcls}\\
            &\ubar{\bm v}_i \leq | V_i | \leq \bar{\bm v}_i ,\quad \forall i \in N, \label{eq:exvbounds}\\
		    &\ubar{\bm \theta}_{ij} \mathcal{I}(V_iV^*_j) \leq \mathcal{R}(V_iV^*_j) \leq \bar{\bm \theta}_{ij} \mathcal{I}(V_iV^*_j) ,\quad \forall (i,j) \in E, \label{eq:exanglediff} \\
            &\ubar{\bm S}_i^g \leq S^g_i \leq \bar{\bm S}_i^g ,\quad \forall i \in N, \label{eq:exsbounds} \\
            &| S_{ij} |^2 \leq \bm t_{ij} ,\quad \forall (i,j), (j,i) \in E. \label{eq:exthermal}
		\end{align}
	\end{subequations}
	}
\end{model}
\section{Convex Relaxations}\label{sec:SDP}
\subsection{The Lasserre Hierarchy} Model \ref{model:ext_ac_opf} is a particular instance of a polynomial optimization problem. In \cite{Lasserre2001}, Lasserre proposed a general purpose approach to approximate the optimal value of such problems via solving a sequence of SDPs of increasing size. The SDP relaxation of quadratic optimization problems (Shor relaxation \cite{Shor1987}, or \cite{Low} in the context of ACOPF) can be understood as the first level of this hierarchy. On higher levels, more lifted variables are introduced and arranged in matrices which are constrained to be positive semi-definite. Lasserre proved that under mild assumptions, this hierarchy of SDP relaxations converges towards the optimal value of the polynomial optimization problem. However, going higher in the hierarchy requires solving SDPs with increasingly large matrix variables. As a consequence, the computational cost becomes prohibitive when the number of polynomial variables increases or when the level of the hierarchy is large. Techniques using chordal sparsity \cite{Waki06sparse} can help extend the reach of such methods, especially for higher levels of the hierarchy. 
In the ACOPF context, \cite{JMPG2015_opf_hierarchy} have shown that the second level Lasserre Hierarchy succeeds in finding the global optimum of the ACOPF in many cases. 

\subsection{The Determinant Hierarchy}
\begin{figure}[h]
    \centering
    \includegraphics[width=0.50\textwidth]{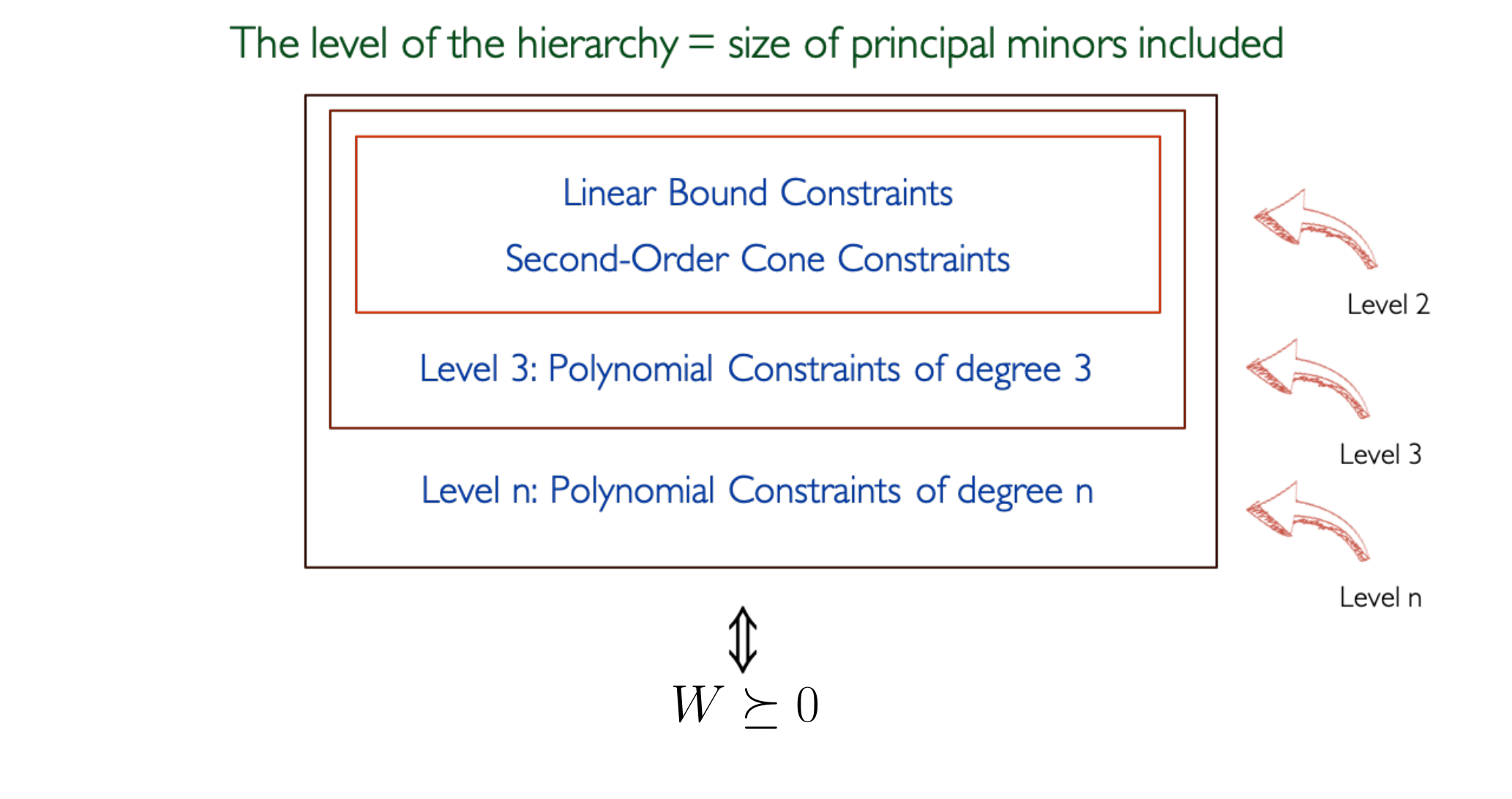}
    \caption{The Determinant Hierarchy.}
    \label{fig:det_hierarchy}
\end{figure}

For a given level of the Lasserre Hierarchy, based on the determinant constraints introduced in \cite{polysdp}, one can define a series of convex relaxations that converge to the corresponding hierarchy level. These relaxations represent the different levels of the Determinant Hierarchy. For a given submatrix of size $k$, one can derive a set of determinant constraints defining level $k$ of the Determinant Hierarchy as depicted in Figure \ref{fig:det_hierarchy}.

In the Determinant Hierarchy, one can trade-off computational efficiency and relaxation strength. In our numerical experiments, we show that the Level 3 Determinant Hierarchy is sufficient to capture the strength of the Level 1 Lasserre Hierarchy on all tested benchmarks. Note that the Level 2 Determinant Hierarchy corresponds to the SOCP relaxation introduced by Jabr \cite{Jabr}.

In order to formulate the various convex relaxations, we introduce the Hermitian matrix $W$ where each entry $W_{ij}$ corresponds to the product $V_iV_j^{*}$ appearing in the non-convex Model \ref{model:ext_ac_opf}.
The Level n Determinant Hierarchy is presented in Model \ref{model:sdp_opf}:
\begin{model}[!h]
	\caption{Det-SDP}
	\label{model:sdp_opf}
	{\scriptsize
	\begin{subequations}
		\vspace{-0.2cm}
		\begin{align}
		& \mbox{\bf minimize: } \nonumber \\ 
		& \sum\limits_{i \in N} \bm c_{0i} + \bm c_{1i}\mathcal{R}(S^g_i) + \left(\bm c_{2i}\mathcal{R}(S^g_i)\right)^2 \nonumber\\
		& \mbox{\bf subject to: } \nonumber  \\
		& S_{ij} = (\bm Y^*_{ij}+ \bm Y^{\bm c^*}_{ij})\frac{W_{ii}}{| \bm T_{ij}|^2}-\bm Y^*_{ij}\frac{W_{ij}}{\bm T_{ij}}  ,\quad \forall (i,j) \in E \label{eq:sdp_pf_from} \\
		& S_{ji} = (\bm Y^*_{ij}+\bm Y^{c^*}_{ij})W_{jj}-\bm Y^*_{ij}\frac{W_{ij}^*}{\bm T^*_{ij}} ,\quad \forall (i,j) \in E \label{eq:sdp_pf_to}\\
        & S_i^g - \bm{S}_i^d -\bm{Y^s}_i W_{ii} = \displaystyle\sum_{(i,j),  (j,i) \in E} S_{ij} ,\quad \forall i \in N, \label{eq:sdpkcls}\\
		& \mathcal{I}(W_{ij})\tan \ubar{\bm \theta}_{ij} \le \mathcal{R}(W_{ij}) \le \mathcal{I}(W_{ij}) \tan \bar{\bm \theta}_{ij},\quad \forall{(i,j) \in {E}}, \label{eq:sdp_angle}\\
        &\ubar{\bm v}^2_i \leq |W_{ii}| \leq \bar{\bm v}^2_i ,\quad \forall i \in N, \label{eq:sdpvbounds}\\
        &\bar{\bm v}_i\bar{\bm v}_j \cos \ubar{\bm \theta}_{ij} \leq \mathcal{R}(W_{ij}) \leq \bar{\bm v}_i\bar{\bm v}_j ,\quad \forall (i,j) \in E, \label{eq:sdpRwijbounds}\\
        &-\bar{\bm v}_i\bar{\bm v}_j \leq \mathcal{I}(W_{ij}) \leq \bar{\bm v}_i\bar{\bm v}_j ,\quad \forall (i,j) \in E, \label{eq:sdpIwijbounds}\\
        & \textrm{det}(W_s) \geq 0 \quad \forall W_s \subseteq W_c \in T({\mathcal G}),  \label{eq:psd_poly}\\
        &\eqref{eq:exsbounds},\eqref{eq:exthermal}.\notag
		\end{align}
	\end{subequations}
	}
\end{model}

where $\mathbf{conv}(f(x)=0)$ represents the automatic convexification of $f(x)=0$ implemented in \href{https://github.com/coin-or/Gravity}{Gravity} \cite{Gravity}, based on expanding the complex products, splitting the constraints into real and imaginary parts, adding convex sides of the equation to the model and using McCormick envelopes \cite{McCormick1976} for convexifying bilinear terms. 

Given a tree decomposition $T({\mathcal G})$ of ${\mathcal G}$, $W_c$ is the submatrix of $W$ formed by selecting the rows and columns corresponding to nodes appearing in clique $C$ of $T({\mathcal G})$.
$W_s$ represent the exponentially many principal minors corresponding to submatrices of $W_c$ (details can be found in \cite{polysdp}).

\subsection{Using the Reformulation Linearization Technique (RLT)}
 The inequalities developed in this work build on previous efforts in \cite{current_lit}. Several studies, e.g.,  \cite{current_lit, qc_mpc_17}, use relationships between the current magnitude and the power flows to strengthen existing ACOPF relaxations. We present a strengthened version of these valid cuts based on the Reformulation Linearization Technique (RLT) \cite{RLT} (see last paragraph). 

Consider the current $I_{ij} \in \mathcal{C}$ which flows from bus $i$ into edge $(i,j)$. 
$I_{ij}$ can be expressed as a linear function of the difference in voltages between adjacent buses and can be used to describe power flows:
 \begin{subequations}
 \begin{align}
& I_{ij}=(\bm Y_{ij}+\bm{Y^c}_{ij})V_i/\bm{T}_{ij}-\bm Y_{ij}V_j ,\quad \forall (i,j) \in E \label{eq:current_ext}\\
& S_{ij} =V_{i}I_{ij}^*/\bm T_{ij} ,\quad \forall (i,j) \in E \label{eq:power_ext}
 \end{align}
 \end{subequations}

Using the RLT approach, we multiply constraints \eqref{eq:current_ext} and \eqref{eq:power_ext} with their conjugate to obtain the following equations:
\begin{subequations}
\begin{align}
 &I_{ij}I_{ij}^*| \bm T_{ij}| ^2=|(\bm Y_{ij}+\bm{Y^c}_{ij})|^2 V_{i}V_i^* -\bm{Y}^*_{ij}(\bm{Y}_{ij}+\bm{Y^c}_{ij})\bm{T}^*_{ij}V_{i}V_j^*\nonumber \\
&\;-\bm{Y}_{ij}(\bm{Y}_{ij} +\bm{Y^c}_{ij})^*\bm{T}_{ij}V_{i}^*V_j +| \bm{Y}_{ij} \bm{T}_{ij}| ^2V_{j}V_j^*,\label{eq:ext_lij_nc}\\
 &| \bm T_{ij}|^2 S_{ij}S_{ij}^* = {I}_{ij}{I}_{ij}^*V_iV_i^*, \quad \forall (i,j) \in E \label{eq:ext_lWij_nc}
\end{align}
\end{subequations}
which we then lift to get the following linear RLT constraints:
\begin{subequations}
\begin{align}
 &l_{ij}| \bm T_{ij}| ^2=|(\bm Y_{ij}+\bm{Y^c}_{ij})|^2 W_{ii} -\bm{Y}^*_{ij}(\bm{Y}_{ij}+\bm{Y^c}_{ij})\bm{T}^*_{ij}W_{ij}\nonumber \\
&\;-\bm{Y}_{ij}(\bm{Y}_{ij} +\bm{Y^c}_{ij})^*\bm{T}_{ij}W_{ij}^*+| \bm{Y}_{ij} \bm{T}_{ij}| ^2W_{jj},\quad \forall (i,j) \in E \label{eq:ext_lij}\\
 &\widehat{lW}_{ij} =  | \bm T_{ij}|^2 (\widehat{p}_{ij}+\widehat{q}_{ij}), \quad \forall (i,j) \in E \label{eq:ext_lWij}
\end{align}
\end{subequations}
where lifted variables $\widehat{p}_{ij}$, $\widehat{q}_{ij}$, and $\widehat{lW}_{ij}$ represent quadratic terms ${p}^2_{ij}$, ${q}^2_{ij}$ and product $l_{ij}W_{ii}$, respectively.
Similarly, we use the expressions for $I_{ji}$ and $S_{ji}$:
 \begin{subequations}
 \begin{align}
& I_{ji}=(\bm Y_{ij}+\bm{Y^c}_{ij})V_j-\bm Y_{ij}V_j/\bm{T}_{ij} ,\quad \forall (i,j) \in E \label{eq:ext_Iji}\\
& S_{ji} =V_{j}I_{ji}^* ,\quad \forall (i,j) \in E\label{eq:ext_Sji}
 \end{align}
 \end{subequations}
 to derive the following RLT cuts:
\begin{subequations}
\begin{align}
 l_{ji}| \bm T_{ij}|^2&=|(\bm Y_{ij}+\bm{Y^c}_{ij}) \bm T_{ij}|^2 W_{jj}-\bm Y_{ij}(\bm Y_{ij}+\bm{Y^c}_{ij})^*\bm T_{ij}^*W_{ij} \nonumber \\
&-\bm Y^*_{ij}(\bm Y_{ij}+\bm{Y^c}_{ij})\bm{T}_{ij}W_{ij}^*+| \bm{Y}_{ij}|^2 W_{ii} \forall (i,j) \in E, \label{eq:ext_lji}\\
  \widehat{lW}_{ji} &= \widehat{p}_{ji}+\widehat{q}_{ji} ,\quad \forall (i,j) \in E, \label{eq:ext_lWji}
\end{align}
\end{subequations}
Note that products we lifted in \eqref{eq:ext_lij}-\eqref{eq:ext_lWji} also appear in other constraints (e.g., the thermal limits constraints).
Indeed, \eqref{eq:exthermal} may be reformulated using the lifted variables $\widehat{p}_{ij}$ and $\widehat{q}_{ij}$ appearing in the RLT cuts to yield:
\begin{align*}
& \widehat{p}_{ij}+\widehat{q}_{ij} \le \bm t_{ij} \quad \forall (i,j);(j,i) \in E
\end{align*}
This linking strengthens the RLT cuts and helps reduce the gap at the root node (pre-bound tightening) on a number of open instances.

The convex relaxation for the extended ACOPF problem is given in Model \ref{model:ext_sdp_current}:

\begin{model}[!h]
	\caption{Det-SDP with RLT}
	\label{model:ext_sdp_current}
	{\scriptsize
	\begin{subequations}
		\vspace{-0.2cm}
		\begin{align}
		& \mbox{\bf minimize: } \nonumber \\ 
		& \sum\limits_{i \in N} \bm c_{0i} + \bm c_{1i}\mathcal{R}(S^g_i) + \left(\bm c_{2i}\mathcal{R}(S^g_i)\right)^2\\
		& \mbox{\bf subject to: } \nonumber  \\
		& \widehat{p}_{ij}+\widehat{q}_{ij} \le \bm t_{ij} \quad \forall (i,j);(j,i) \in E \label{eq:lifted_therm}\\
		& {\mathcal Mc}\left(\widehat{Wl}_{ij},W_{ii}, l_{ij}\right) ,\quad  \forall (i,j);(j,i) \in E\label{eq:Mc}\\
		& {\mathcal Sec}\left(\widehat{p}_{ij},p_{ij}\right),\;{\mathcal Sec}\left(\widehat{q}_{ij},q_{ij}\right) ,\quad  \forall (i,j);(j,i) \in E\label{eq:Sec}\\
		&\eqref{eq:exsbounds},\eqref{eq:sdp_pf_from}-\eqref{eq:psd_poly},\eqref{eq:ext_lij}-\eqref{eq:ext_lWij},\eqref{eq:ext_lji}-\eqref{eq:ext_lWji},\eqref{lub}.\notag
		\end{align}
	\end{subequations}
	}
\end{model}

where ${\mathcal Sec}(y, x)$ denotes the convex envelope of the set $S=\left\{(y,x) \in \mathbf{R}^2: y=x^2, ~\ubar{\bm x} \le x \le  \bar{\bm x}\right\}$ defined by the secant inequality $y-x(\ubar{\bm x}+\bar{\bm x})+\ubar{\bm x}\bar{\bm x} \le 0$, and the quadratic constraint $y \ge x^2$.
Note that the quality of the McCormick relaxation ${\mathcal Mc}\left(\widehat{Wl}_{ij},W_{ij}, l_{ij}\right)$ strongly depends on the tightness of the lower and upper bounds on the original variables. As suggested in \cite{qc_mpc_17} we can use the following initial bounds on $l_{ij}$ and $l_{ji}$:
\begin{align}
    &  0 \le l_{ij} \le \bm{t}_{ij}/\ubar{\bm v}^2_i \quad \forall (i,j);(j,i) \in E \label{lub}
\end{align}

\section{Semidefinite Programming-based Bound Tightening (SDP-BT)} \label{sec:SDP-BT}

\subsection{Optimality-Based Bound Tightening (OBBT)}\label{sec:obbt}
OBBT \cite{Zamora, Gleixner} is a fundamental technique used in global optimization algorithms. 
The algorithm is used to tighten lower and upper bounds for each variable $x_i$ that appears in a given convex relaxation of the original problem. On the one hand, tighter bounds reduce the domain to be explored for global optimization. On the other hand, the quality of convex relaxations critically depends on the tightness of the bounds. This can lead to a positive cycle where the bounds get smaller which leads to tighter relaxations which in turn can reduce the bounds, until a fixed-point is reached. Details about the properties of this algorithm can be found in \cite{CP2015}.
The domain of variable $x_i$ may be reduced by solving a minimization (resp. maximization) OBBT-subproblem which identifies $\ubar{x}_i$ (resp. $\bar{x}_i$), the lowest (resp. highest) value of $x_i$ such that the objective function $f$ is upper bounded by $\bar{f}$, the objective value at a given feasible solution. 
OBBT algorithms tailored to the ACOPF have been developed in \cite{acopf_obbt}. More recently, \cite{lu2018tight,Bynum19, KDS2018minors} looked at the QC-based \cite{qc_mpc_17}-\cite{qc} OBBT.
\subsection{SDP-based Bound Tightening (SDP-BT)}
The pair of subproblems corresponding to a variable $x_i$ in Model \ref{model:ext_sdp_current} is shown in Model \ref{model:sdp_current_obbt}

\begin{model}[ht]
	\caption{SDP-BT subproblem}
	\label{model:sdp_current_obbt}
	{\scriptsize
	\begin{subequations}
		\vspace{-0.2cm}
		\begin{align}
		& \mbox{\bf min/max: } \nonumber \\ 
		&  x_i\nonumber\\
		& \mbox{\bf subject to: } \nonumber  \\
		& \sum\limits_{i \in N} \bm c_{0i} + \bm c_{1i}\mathcal{R}(S^g_i) + \left(\bm c_{2i}\mathcal{R}(S^g_i)\right)^2 \le \bar{f},\label{eq:obj_ub}\\
		&\eqref{eq:exsbounds},\eqref{eq:sdp_pf_from}-\eqref{eq:psd_poly},\eqref{eq:ext_lij}-\eqref{eq:ext_lWij},\eqref{eq:ext_lji}-\eqref{eq:ext_lWji},\eqref{eq:lifted_therm}-\eqref{eq:Sec},\eqref{lub}.\notag
		\end{align}
	\end{subequations}
	}
\end{model}

Each SDP-BT subproblem (shown in Model \ref{model:sdp_current_obbt}) is nonlinear but defines a convex feasible region \cite{polysdp}. Note that unlike previous approaches (such as \cite{obbt_Coffrin}), we use a tight relaxation in each bound tightening subproblem using the Determinant Hierarchy discussed previously. 

The determinant-cuts representation can leverage scalable nonlinear optimization tools such as Ipopt \cite{Ipopt}, leading to tight and computationally efficient models to be used in SDP-BT.

In each major iteration $k$ of SDP-BT, minimization and maximization sub-problems $M_i$ (based on Model \ref{model:sdp_current_obbt}) are solved for a variable $x_i$ until its domain satisfies conditions \eqref{C1}-\eqref{C2}:
\begin{gather}
\left(\left|\bar {x}_i^{(k)}-\bar {x}_i^{(k-1)}\right| \le {\epsilon_d} \textrm{ AND } \left| \ubar {x}_i^{(k)}-\ubar {x}_i^{(k-1)}\right| \le {\epsilon_d} \right)\label{C1} \\ \textrm{OR} \nonumber \\\left| \bar {x}_i^{(k)}-\ubar {x}_i^{(k)} \right| \le {\epsilon_d} \label{C2},
\end{gather}
where $\epsilon_d >0$.
Variables are grouped into batches that run on separate machines and corresponding bounds are dynamically updated.

The overall algorithm is presented in Algorithm \ref{obbt_algo}, where $gap$ is defined as $(\bar{f}-\ubar{f})/\bar{f}$.

\begin{algorithm}{} 
\begin{algorithmic}[1]
\State Set SDP-BT convergence tolerances $\epsilon_d$ and $\epsilon_o$, maximum number of iterations $\bm {K}$ and total number of threads $\bm{T}$.
\State $k$ = 0, $end$ = false, $\mathcal{M} = \emptyset$ and $x_i(close)$ = \textrm{false}, $\forall i \in \{1,\ldots, n\}$, $gap = 100$.
\State Solve initial relaxation, get $\ubar{f}$ and update $gap$ \label{step_gap}
\While {($gap \ge \epsilon_o \textrm{ OR } end = \textrm{false}) \textrm{ AND } k \le \bm {K}$}
\State $end$ = true
\For {$i \in \{1,\ldots, n\}$}
    \If {$x_i(close)$ = \textrm{false}}
        \State Add $M_i$ to $\mathcal{M}$
        \If {$\textrm{size}(\mathcal{M})=\bm{T}$}
            \State Run all $\mathcal{M}$ models in parallel
            \State update variable bounds on ${x}_j \in M_j$
            \For {$M_j \in \mathcal{M}$}
                \If{\eqref{C1}-\eqref{C2} holds for $x_j$}
                \State $x_j(close)$ = true
                \Else
                     \State $end$ = false
                     \EndIf
                        \EndFor
            \State $\mathcal{M} = \emptyset$
        \EndIf
        \EndIf
\EndFor
\State $k=k+1$, Compute new $\ubar{f}$ and update $gap$
\EndWhile
\end{algorithmic}
\caption{SDP-BT}
\label{obbt_algo}
\end{algorithm}

\section{Numerical experiments}
\label{sec:numeric}
In this section we test the convex relaxations developed in this work on case studies from benchmark libraries \cite{pglib} and \cite{nesta}. First, in section \ref{sec:casestudy9}, we consider a small test case from \cite{nesta} that exhibits a large optimality gap and highlights the importance of combining RLT cuts with the Determinant Hierarchy.

In section \ref{sec:pglib}, we show the results of our approach on all problems having up to 300 nodes in PGLIB 19.05 \cite{pglib}.

Finally, in section \ref{sec:nesta3}, we compare the performance of SDP-BT with the gaps reported in \cite{KDS2018minors} on NESTA 0.3.

For the first level of the Lasserre Hierarchy we use the Julia package PowerModels.jl \cite{powermodels} which implements the standard SDP relaxation strengthened by lifted non-linear cuts (LNC) introduced in \cite{LNC}. For the second level we use our own preliminary Julia implementation which expands SumOfSquares.jl \cite{SumOfSquares} to use chordal sparsity. In both cases, we use version 8 of Mosek \cite{mosek} with standard parameter setting to solve the resulting SDPs.

Our implementation of the Lasserre Hierarchy is a general purpose tool, i.e., it is not tailored for the ACOPF problem. In particular, we do not increase the level selectively as done in \cite{MH2015_selective_hierarchy}. We also emphasize that we do not change the objective function, nor the input data (e.g., changing values of low impedance lines) when formulating the relaxations. 
Finally, we note that our implementation of the Lasserre Hierarchy uses the following lifting: $\tilde{V}\tilde{V}^T$ where $\tilde{V}_\alpha = (\mathcal{R}V_{\alpha_1}\mathcal{I}V_{\alpha_2})$ for $\alpha \in \{(i,j) \mid i+j \leq 2\}$. We refer to \cite{CM2018_complex_hierarchy} for a detailed discussion on the comparison to the lifting used for the other implementations in this paper. 

 Algorithm\,\ref{obbt_algo} is implemented in C++ within Gravity \cite{Gravity} and linked to Intel's Message Passing Interface (MPI) library 19.0.4. All optimization problems are solved using Ipopt \cite{Ipopt}, with the linear solver Ma27 \cite{HSL}, using an absolute constraint satisfaction tolerance of 1e-6. To solve this parallel implementation, we use 12 machines with 12 threads each. The maximum time per instance is set to 6 hours. The optimality tolerance of the SDP-BT algorithm, $\bm{\epsilon_o}$ is set to 1\% and the domain tolerance $\bm{\epsilon_d}$ is set to 1e-3. 
 The code is open-source and can be accessed at \url{https://github.com/coin-or/Gravity/tree/master/examples/Optimization/NonLinear/Power}.

\subsection{Case study: case 9}\label{sec:casestudy9}
\begin{figure}
    \centering
    \includegraphics[width=0.43\textwidth]{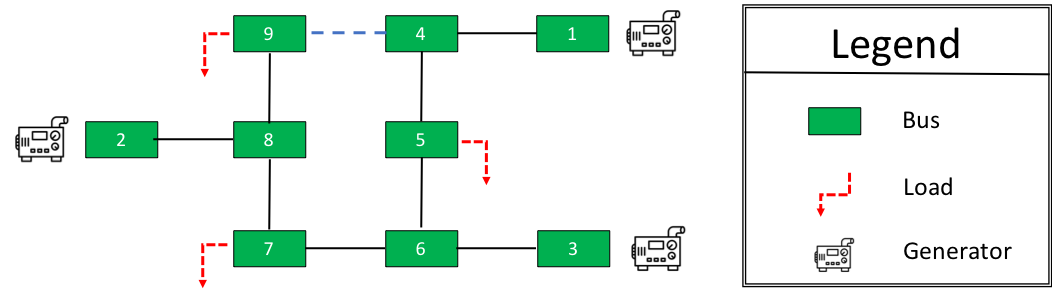}
    \caption{Case 9}
    \label{fig:case9}
\end{figure}

In this section, we discuss nesta\_case9\_bgm\_\_nco \cite{Bukhsh} referred to as case9, and depicted in \figureautorefname{ \ref{fig:case9}}. 
This instance has 4 locally optimal solutions \cite{Bukhsh} and an SDP (Lasserre Hierarchy Level 1) optimality gap of 10.84\% \cite{polysdp}. 

The underlying graph contains only one cycle and thus disconnecting any edge of the cycle yields a tree. 
For the tree produced by omitting edge $(4,9)$ (dotted edge in \figureautorefname{ \ref{fig:case9}}), the root-node optimality gap using the first level of the Lasserre Hierarchy is 5.33 \%.  We refer to this instance as case9\_tree. 

We first study the effect of the SDP determinant constraints and the valid RLT cuts on the SDP-BT algorithm. We compare a version of SDP-BT excluding the RLT cuts, a version using the SOCP relaxation \cite{Jabr} including the RLT cuts and the proposed version based on Model \ref{model:ext_sdp_current}. These are also compared to the Lasserre Hierarchy Levels 1 and 2. Results are summarized in Table~\ref{tab:case9}. The large gaps obtained by adding only the RLT cuts or only the SDP cuts highlight the importance of combining both for strengthening convex relaxations of the ACOPF problem. 
\subsection{Why combine RLT with SDP cuts?}
Let us consider the simple case of the RLT equation \eqref{eq:ext_lij} in the absence of line charging and transformer parameters.
\begin{align}
l_{ij} =(g_{ij}^2+b^2_{ij}) (W_{ii}+W_{jj}-2\mathcal{R}(W_{ij}))
\end{align}
Note that this is only a scaled version of the real power loss on an edge,
\begin{align}
p_{ij}+p_{ji} =g_{ij} (W_{ii}+W_{jj}-2\mathcal{R}(W_{ij}))
\end{align}
Thus, an upper bound on $l_{ij}$ implies an upper bound on the power loss. 
Conversely, the determinant constraints provide a lower bound on the loss. This is easy to see at level 2 of the Determinant Hierarchy with the second order cone constraint,
\begin{align}
W_{ii}W_{jj} \ge (\mathcal{R}(W)_{ij})^2+(\mathcal{I}(W)_{ij})^2.
\end{align}
Using the inequality of arithmetic and geometric means:
\begin{align}
(W_{ii}+W_{jj})^2/4 \ge W_{ii}W_{jj} \ge \mathcal{R}(W)_{ij}^2\\
\implies W_{ii}+W_{jj} - 2\mathcal{R}(W_{ij}) \ge 0.
\end{align}
We conjecture that higher order determinant constraints further tighten the lower bounds on the power loss and that performing bound tightening triggers a positive cycle where the relaxation gets stronger and the bounds get smaller in an alternating fashion.

\begin{table}
\centering
\begin{tabular}{|c|rr|rr|}
\hline 
 &    \multicolumn{2}{c|}{case9}   & \multicolumn{2}{c|}{case9\_tree} \\
 & gap & time & gap &time\\
\hline 
OBBT \& RLT &  10.20 & 10.27 & 0.66 & 6.00\\
OBBT \& SDP no RLT & 10.82 & 8.07 & 5.33 & 6.89\\
OBBT \& SDP \& RLT (SDP-BT) & \bf 0.44 & \bf 10.19 & \bf 0.66 & \bf 6.03\\
\hline
\end{tabular}
\vspace{0.25cm}
\caption{Optimality gaps and computation time for different relaxations on case9 and case9\_tree.\label{tab:case9}}
\end{table}


\subsection{PGLIB 19.05} \label{sec:pglib}

We run Algorithm\,\ref{obbt_algo} on all instances in \cite{pglib} having less than 300 buses.
The results from the numerical experiments are shown in Table \ref{tab: results_pglib} where instances with an initial gap $>1\%$ are highlighted in {\bf bold}. For each case we report the root gap before calling SDP-BT, the final gap after running SDP-BT, as well as the number of iterations and the total wall clock time (in seconds) of the algorithm (T.L. indicating that we hit the time limit).
We also report the level 1 and 2 Lasserre Hierarchy gaps (gap1 and gap2) as well as the corresponding total wall-clock time. Note that level 1 Lasserre Hierarchy relaxation includes the lifted nonlinear cuts (LNC) introduced in \cite{LNC}. Missing entries correspond to instances where the code ran out of memory. Starred entries refer to Mosek reporting unknown primal and dual status. 
The SDP-BT algorithm manages to decrease the optimality gap to less than 1\% {\bf on all tested instances}. 
The root node gap from Model \ref{model:ext_sdp_current} is tight and closes the gap in 38 out of 48 instances. 
It is  worth mentioning that the number of iterations taken by the SDP-BT algorithm is typically small, taking only one iteration to converge on all but 2 instances. 

\subsection{NESTA 0.3 \cite{KDS2018minors}}
\label{sec:nesta3}
In this section, we present an extensive comparison of SDP-BT's performance, on all instances having less than 300 nodes with state-of-the-art results published in \cite{KDS2018minors}.
The results are reported in Table \ref{tab:nesta3}. In the last four columns, we report optimality gaps and computation times for both the root node relaxation (including OBBT + cut generation) and the Branch\&Cut algorithm developed in \cite{KDS2018minors}. The SDP-BT algorithm closes the gap on 55 out of 57 instances, outperforming the root node relaxation reported in \cite{KDS2018minors} on all instances and closing 4 instances more than the Branch\&Cut algorithm.

\section{Conclusion and Future Work}

Parallelizing global optimization algorithms is not a trivial task. The SDP-BT algorithm and the underlying code we present here relies solely on open-source tools and can take full advantage of computational clusters using state-of-the-art MPI libraries. In our experiments, we showcase the performance of SDP-BT on a small cluster of 10 machines, closing the gap on a number of open test cases from the literature.


Several avenues for further investigation remain. The level 2 Lasserre Hierarchy, which involves the solution of large SDPs, may be solved efficiently using the Determinant Hierarchy presented here. 
We are also looking at generating outer-approximation cuts that would lead to linear Determinant Hierarchies. Finally, we are working on integrating our contributions into an adaptive variable partitioning algorithm \cite{amp} which will allow us to tackle larger instances.
{
\begin{table}	
\resizebox{0.5\textwidth}{!}{%
\begin{tabular}{|l|rrrr|rrrr|}																	
																	
\hline																	
	&	\multicolumn{4}{c|}{SDP-BT}							&	\multicolumn{4}{c|}{Lasserre relaxation}							\\
Case	&	root gap	&	final gap	&	time	&	iter	&	gap1	&	time1	&	gap2	&	time2	\\
\hline																	
				
pglib\_opf\_case3\_lmbd	&	0.00	&	0.00	&	0.08	&	0	&	0.38	&	0.01	&	0.00	&	0.23	\\
\bf pglib\_opf\_case5\_pjm	&	\bf 0.09	& \bf	0.09	& \bf	0.13	& \bf	0	& \bf	5.22	& \bf	0.01	& \bf	0.00	& \bf	20.87	\\
pglib\_opf\_case14\_ieee	&	0.00	&	0.00	&	0.50	&	0	&	0.00	&	0.22	&	0.00	&	141.99	\\
pglib\_opf\_case24\_ieee\_rts	&	0.00	&	0.00	&	0.86	&	0	&	0.00	&	0.19	&	-	&	-	\\
pglib\_opf\_case30\_as	&	0.00	&	0.00	&	0.85	&	0	&	0.00	&	0.14	&	0.00	&	5924.40	\\
pglib\_opf\_case30\_fsr	&	0.00	&	0.00	&	0.71	&	0	&	0.01	&	0.23	&	-	&	-	\\
pglib\_opf\_case30\_ieee	&	0.00	&	0.00	&	0.68	&	0	&	0.02	&	0.17	&	-	&	-	\\
pglib\_opf\_case39\_epri	&	0.00	&	0.00	&	1.47	&	0	&	0.01	&	0.23	&	0.00	&	1336.31	\\
pglib\_opf\_case57\_ieee	&	0.00	&	0.00	&	1.62	&	0	&	0.01	&	0.40	&	-	&	-	\\
pglib\_opf\_case73\_ieee\_rts	&	0.00	&	0.00	&	3.54	&	0	&	0.00	&	0.60	&	-	&	-	\\
pglib\_opf\_case89\_pegase	&	0.29	&	0.29	&	12.96	&	0	&	0.34	&	1.97	&	-	&	-	\\
pglib\_opf\_case118\_ieee	&	0.03	&	0.03	&	3.18	&	0	&	0.07	&	1.29	&	-	&	-	\\
\bf pglib\_opf\_case162\_ieee\_dtc	&	\bf 1.57	&	\bf 0.45	&	\bf 5723.00	&	\bf 1	&	\bf 1.78	&	\bf 5.83	&	\bf -	&	\bf -	\\
pglib\_opf\_case179\_goc	&	0.07	&	0.07	&	6.98	&	0	&	0.07	&	1.64	&	-	&	-	\\
pglib\_opf\_case200\_tamu	&	0.00	&	0.00	&	7.11	&	0	&	0.00	&	1.47	&	-	&	-	\\
\bf pglib\_opf\_case300\_ieee	&	\bf 0.10	&	\bf 0.10	&	\bf 18.99	&	\bf 0	&	\bf 1.56	&	\bf 3.45	&	\bf -	&	\bf -	\\
\hline																	
																	
\bf pglib\_opf\_case3\_lmbd\_\_api	&	\bf 0.93	&	\bf 0.93	&	\bf 0.10	&	\bf 0	&	\bf 4.99	&	\bf 0.01	&	\bf 0.00	&	\bf 0.20	\\
pglib\_opf\_case5\_pjm\_\_api	&	0.01	&	0.01	&	0.20	&	0	&	0.30$^*$	&	0.04$^*$	&	0.00	&	19.33	\\
pglib\_opf\_case14\_ieee\_\_api	&	0.01	&	0.01	&	0.62	&	0	&	0.02	&	0.06	&	0.00	&	143.71	\\
\bf pglib\_opf\_case24\_ieee\_rts\_\_api	&	\bf 1.03	&	\bf 0.03	&	\bf 21.42	&	\bf 1	&	\bf 2.07	&	\bf 0.21	&	\bf -	&	\bf -	\\
pglib\_opf\_case30\_as\_\_api	&	0.72	&	0.72	&	0.82	&	0	&	16.19$^*$	&	0.21$^*$	&	0.39	&	9807.31	\\
pglib\_opf\_case30\_fsr\_\_api	&	0.27	&	0.27	&	2.24	&	0	&	0.52	&	0.20	&	-	&	-	\\
pglib\_opf\_case30\_ieee\_\_api	&	0.02	&	0.02	&	0.74	&	0	&	0.34$^*$	&	0.17$^*$	&	-	&	-	\\
pglib\_opf\_case39\_epri\_\_api	&	0.16	&	0.16	&	0.71	&	0	&	0.46	&	0.26	&	0.01	&	1973.40	\\
pglib\_opf\_case57\_ieee\_\_api	&	0.00	&	0.00	&	3.34	&	0	&	0.02	&	0.51	&	-	&	-	\\
\bf pglib\_opf\_case73\_ieee\_rts\_\_api	&	\bf 2.13	&	\bf 0.75	&	\bf 261.14	&	\bf 1	&	\bf 2.92	&	\bf 0.65	&	\bf -	&	\bf -	\\
\bf pglib\_opf\_case89\_pegase\_\_api	&	\bf 11.70	&	\bf 0.93	&\bf 	6013.07	&	\bf 3	&	\bf 12.12	&	\bf 2.27	&	\bf -	&	\bf -	\\
\bf pglib\_opf\_case118\_ieee\_\_api	&	\bf 8.44	&	\bf 0.99	&	\bf 2030.47	&	\bf 5	&	\bf 11.20	&	\bf 1.26	&\bf 	-	&	\bf -	\\
\bf pglib\_opf\_case162\_ieee\_dtc\_\_api	&	\bf 1.26	&	\bf 0.26	&	\bf 16277.76	&	\bf 1	&	\bf 1.44	&	\bf 5.33	&	\bf -	&	\bf -	\\
pglib\_opf\_case179\_goc\_\_api	&	0.54	&	0.54	&	9.01	&	0	&	0.55	&	1.35	&	-	&	-	\\
pglib\_opf\_case200\_tamu\_\_api	&	0.00	&	0.00	&	39.17	&	0	&	0.00	&	2.05	&	-	&	-	\\
pglib\_opf\_case300\_ieee\_\_api	&	0.07	&	0.07	&	21.87	&	0	&	0.21	&	3.45	&	-	&	-	\\
\hline																	
																	
pglib\_opf\_case3\_lmbd\_\_sad	&	0.10	&	0.10	&	0.08	&	0	&	0.62	&	0.01	&	0.00	&	0.21	\\
pglib\_opf\_case5\_pjm\_\_sad	&	0.00	&	0.00	&	0.20	&	0	&	0.00	&	0.03	&	0.00	&	18.36	\\
pglib\_opf\_case14\_ieee\_\_sad	&	0.11	&	0.11	&	0.36	&	0	&	0.09	&	0.10	&	0.00	&	148.34	\\
\bf pglib\_opf\_case24\_ieee\_rts\_\_sad	&	\bf 3.54	&	\bf 0.11	&	\bf 21.69	&	\bf 1	&	\bf 2.52	&	\bf 0.17	&	\bf -	&	\bf -	\\
pglib\_opf\_case30\_as\_\_sad	&	0.21	&	0.21	&	0.95	&	0	&	0.16	&	0.25	&	0.00	&	6529.10	\\
pglib\_opf\_case30\_fsr\_\_sad	&	0.02	&	0.02	&	0.69	&	0	&	0.02	&	0.19	&	-	&	-	\\
pglib\_opf\_case30\_ieee\_\_sad	&	0.00	&	0.00	&	0.84	&	0	&	0.00	&	0.16	&	-	&	-	\\
pglib\_opf\_case39\_epri\_\_sad	&	0.02	&	0.02	&	1.36	&	0	&	0.02	&	0.24	&	-	&	-	\\
pglib\_opf\_case57\_ieee\_\_sad	&	0.04	&	0.04	&	4.14	&	0	&	0.04	&	0.63	&	-	&	-	\\
\bf pglib\_opf\_case73\_ieee\_rts\_\_sad	&	\bf 2.13	&	\bf 0.33	&	\bf 228.97	&	\bf 1	&	\bf 1.48	&	\bf 0.58	&	\bf -	&	\bf -	\\
pglib\_opf\_case89\_pegase\_\_sad	&	0.29	&	0.29	&	12.35	&	0	&	0.32	&	1.95	&	-	&	-	\\
\bf pglib\_opf\_case118\_ieee\_\_sad	&	\bf 2.49	&	\bf 0.18	&	\bf 471.40	&	\bf 1	&	\bf 1.83	&	\bf 1.21	&	\bf -	&	\bf -	\\
\bf pglib\_opf\_case162\_ieee\_dtc\_\_sad	&	\bf 1.38	&	\bf 0.28	&	\bf 3666.39	&	\bf 1	&	\bf 1.79	&	\bf 5.63	&	\bf -	&	\bf -	\\
pglib\_opf\_case179\_goc\_\_sad	&	0.94	&	0.94	&	15.01	&	0	&	0.91	&	1.59	&	-	&	-	\\
pglib\_opf\_case200\_tamu\_\_sad	&	0.00	&	0.00	&	8.50	&	0	&	0.00	&	1.64	&	-	&	-	\\
\bf pglib\_opf\_case300\_ieee\_\_sad	&	\bf 0.12	&	\bf 0.12	&	\bf 13.98	&	\bf 0	&	\bf 1.40	&	\bf 3.28	&	\bf -	&	\bf -	\\
\hline																	
\end{tabular}																	
}
\vspace{0.2cm}
\caption{PGLIB 19.05 results\label{tab: results_pglib}}
\end{table}																	
		}

		{
\begin{table}	
\resizebox{0.5\textwidth}{!}{%
\begin{tabular}{|l|rrrr|rrrr|}																	
																	
\hline																	
	&	\multicolumn{4}{c|}{SDP-BT}							&	\multicolumn{4}{c|}{Best of \cite{KDS2018minors}}							\\
Case	&	root gap	&	final gap	&	time	&	iter	&	OBBT+cuts gap	&	time1	&	final gap	&	time2	\\
\hline																	
														
nesta\_case3\_lmbd	&	0.00	&	0.00	&	0.03	&	0	&	0.10	&	0.95	&	0.09	&	0.95	\\
nesta\_case4\_gs	&	0.00	&	0.00	&	0.03	&	0	&	0.00	&	0.03	&		&		\\
\bf nesta\_case5\_pjm	&\bf	0.11	&\bf	0.11	&	\bf0.05	&\bf	0	&\bf	2.11	&	\bf3.26	&	\bf0.10	&	\bf 108.39	\\
nesta\_case6\_c	&	0.00	&	0.00	&	0.03	&	0	&	-	&	-	&	-	&	-	\\
nesta\_case6\_ww	&	0.00	&	0.00	&	0.06	&	0	&	0.01	&	1.08	&		&		\\
nesta\_case9\_wscc	&	0.00	&	0.00	&	0.05	&	0	&	0.00	&	0.09	&		&		\\
nesta\_case14\_ieee	&	0.00	&	0.00	&	0.09	&	0	&	0.00	&	2.70	&		&		\\
nesta\_case24\_ieee\_rts	&	0.00	&	0.00	&	0.16	&	0	&	-	&	-	&	-	&	-	\\
nesta\_case29\_edin	&	0.00	&	0.00	&	0.73	&	0	&	0.01	&	33.99	&		&		\\
nesta\_case30\_as	&	0.00	&	0.00	&	0.15	&	0	&	0.06	&	0.11	&		&		\\
nesta\_case30\_fsr	&	0.01	&	0.01	&	0.17	&	0	&	0.07	&	14.49	&		&		\\
nesta\_case30\_ieee	&	0.02	&	0.02	&	0.16	&	0	&	0.03	&	14.55	&		&		\\
nesta\_case39\_epri	&	0.01	&	0.01	&	0.25	&	0	&	0.05	&	0.25	&		&		\\
nesta\_case57\_ieee	&	0.00	&	0.00	&	0.40	&	0	&	0.06	&	0.22	&		&		\\
nesta\_case73\_ieee\_rts	&	0.00	&	0.00	&	2.32	&	0	&	-	&	-	&	-	&	-	\\
nesta\_case118\_ieee	&	0.02	&	0.02	&	0.97	&	0	&	0.14	&	355.50	&	0.10	&	502.47	\\
\bf nesta\_case162\_ieee\_dtc	& \bf	0.88	& \bf	0.88	& \bf	8.00	& \bf	0	& \bf	1.57	& \bf	948.30	& \bf	1.46	& \bf	1837.44	\\
nesta\_case189\_edin	&	0.05	&	0.05	&	1.08	&	0	&	0.04	&	63.15	&		&		\\
nesta\_case300\_ieee	&	0.07	&	0.07	&	3.41	&	0	&	0.09	&	520.50	&		&		\\
\hline																	
																	
nesta\_case3\_lmbd\_\_api	&	0.32	&	0.32	&	0.02	&	0	&	0.81	&	1.05	&	0.02	&	3.34	\\
nesta\_case4\_gs\_\_api	&	0.02	&	0.02	&	0.03	&	0	&	0.03	&	0.55	&		&		\\
nesta\_case5\_pjm\_\_api	&	0.00	&	0.00	&	0.04	&	0	&	0.05	&	0.81	&		&		\\
nesta\_case6\_c\_\_api	&	0.01	&	0.01	&	0.04	&	0	&	-	&	-	&	-	&	-	\\
nesta\_case6\_ww\_\_api	&	0.02	&	0.02	&	0.06	&	0	&	0.00	&	3.39	&		&		\\
nesta\_case9\_wscc\_\_api	&	0.00	&	0.00	&	0.05	&	0	&	0.00	&	0.06	&		&		\\
nesta\_case14\_ieee\_\_api	&	0.05	&	0.05	&	0.08	&	0	&	0.04	&	13.18	&		&		\\
nesta\_case24\_ieee\_rts\_\_api	&	0.54	&	0.54	&	0.20	&	0	&	-	&	-	&	-	&	-	\\
nesta\_case29\_edin\_\_api	&	0.00	&	0.00	&	2.38	&	0	&	0.04	&	136.83	&		&		\\
nesta\_case30\_as\_\_api	&	0.29	&	0.29	&	0.17	&	0	&	0.09	&	62.09	&		&		\\
\bf nesta\_case30\_fsr\_\_api	& \bf	4.93	& \bf	2.66	& \bf	71.69	& \bf	6	& \bf	5.15	& \bf	90.56	& \bf	0.83	& \bf	1802.18	\\
nesta\_case30\_ieee\_\_api	&	0.10	&	0.10	&	0.19	&	0	&	0.06	&	60.03	&		&		\\
nesta\_case39\_epri\_\_api	&	0.00	&	0.00	&	0.53	&	0	&	0.01	&	26.33	&		&		\\
nesta\_case57\_ieee\_\_api	&	0.09	&	0.09	&	0.43	&	0	&	0.06	&	125.44	&		&		\\
nesta\_case73\_ieee\_rts\_\_api	&	0.35	&	0.35	&	1.10	&	0	&	-	&	-	&	-	&	-	\\
\bf nesta\_case118\_ieee\_\_api	&	\bf 17.50	& \bf	1.91	& \bf	2981.79	& \bf	12	& \bf	7.83	& \bf	911.90	& \bf	7.83 \bf	& \bf	1834.74	\\
\bf nesta\_case162\_ieee\_dtc\_\_api	& \bf	0.84	&	\bf 0.84	& \bf	8.03	& \bf	0	& \bf	1.03 	& \bf	2007.66	&	\bf 1.03	& \bf	2007.68	\\
nesta\_case189\_edin\_\_api	&	0.12	&	0.12	&	1.10	&	0	&	0.91	&	592.86	&	0.12	&	663.19	\\
nesta\_case300\_ieee\_\_api	&	0.00	&	0.00	&	53.02	&	0	&	0.10	&	1048.07	&		&		\\
\hline																	
																	
nesta\_case3\_lmbd\_\_sad	&	0.11	&	0.11	&	0.02	&	0	&	0.09	&	1.29	&	0.03	&	1.29	\\
nesta\_case4\_gs\_\_sad	&	0.05	&	0.05	&	0.03	&	0	&	0.01	&	0.66	&		&		\\
nesta\_case5\_pjm\_\_sad	&	0.04	&	0.04	&	0.04	&	0	&	0.07	&	0.94	&		&		\\
nesta\_case6\_c\_\_sad	&	0.01	&	0.01	&	0.04	&	0	&	-	&	-	&	-	&	-	\\
nesta\_case6\_ww\_\_sad	&	0.00	&	0.00	&	0.08	&	0	&	0.00	&	1.53	&		&		\\
nesta\_case9\_wscc\_\_sad	&	0.03	&	0.03	&	0.06	&	0	&	0.01	&	1.14	&		&		\\
nesta\_case14\_ieee\_\_sad	&	0.00	&	0.00	&	0.09	&	0	&	0.06	&	0.16	&		&		\\
\bf nesta\_case24\_ieee\_rts\_\_sad	& \bf	5.13	&	\bf 0.34	& \bf	17.95	& \bf	1	&	-	&	-	&	-	&	-	\\
\bf nesta\_case29\_edin\_\_sad	&	\bf 23.21	& \bf	0.81	& \bf	215.20	& \bf	2	& \bf	0.70	&	\bf 325.68	& \bf	0.67	& \bf	1837.01	\\
nesta\_case30\_as\_\_sad	&	0.47	&	0.47	&	0.25	&	0	&	0.09	&	38.85	&		&		\\
nesta\_case30\_fsr\_\_sad	&	0.10	&	0.10	&	0.18	&	0	&	0.09	&	26.57	&		&		\\
nesta\_case30\_ieee\_\_sad	&	0.03	&	0.03	&	0.15	&	0	&	0.02	&	26.78	&		&		\\
nesta\_case39\_epri\_\_sad	&	0.05	&	0.05	&	0.28	&	0	&	0.02	&	11.54	&		&		\\
nesta\_case57\_ieee\_\_sad	&	0.04	&	0.04	&	0.42	&	0	&	0.07	&	36.75	&		&		\\
\bf nesta\_case73\_ieee\_rts\_\_sad	& \bf	3.42	&	\bf 0.73	& \bf	202.15	& \bf	1	&	-	&	-	&	-	&	-	\\
\bf nesta\_case118\_ieee\_\_sad	& \bf	5.93	& \bf	0.29	& \bf	1062.24	& \bf	2	& \bf	3.35	& \bf	748.42	& \bf	3.07	& \bf	1804.74	\\
\bf nesta\_case162\_ieee\_dtc\_\_sad	& \bf	3.23	&	\bf 0.29	& \bf	T.L.	& \bf	2	& \bf	3.76	&	\bf 1741.94	&		&		\\
\bf nesta\_case189\_edin\_\_sad	&	\bf 1.23	& \bf	0.25	& \bf	562.08	& \bf	1	& \bf	1.41	& \bf	315.67	& \bf	1.06	& \bf	1814.79	\\
nesta\_case300\_ieee\_\_sad	&	0.09	&	0.09	&	3.57	&	0	&	0.10	&	1226.36	&		&		\\
\hline																															
\end{tabular}		
}
\vspace{0.2cm}
\caption{NESTA 0.3 results\label{tab:nesta3} comparing with \cite{KDS2018minors}}
\end{table}																	
		}															




%
\bibliographystyle{IEEEtran}
\bibliography{ACOPF}

\end{document}